\documentclass[12pt]{iopart}

\usepackage{amsfonts,amsmath,amsthm}    
\usepackage{cite}                
\usepackage[dvips]{graphicx}
\usepackage{booktabs}            
\usepackage{titlesec}            
\titleformat{\section}{\normalfont\Large\bfseries}{\thesection}{1em}{}
\titleformat{\subsection}[runin]{\normalfont\bfseries}{\thesubsection.}{0.5em}{}[.]
\usepackage{algorithmic}
\usepackage{algorithm}
\usepackage{booktabs}

\newcommand{\R}{\mathbb{R}}
\newcommand{\N}{\mathbb{N}}

\newcommand{\embed}{\hookrightarrow}

\newcommand{\norm}[1]{\|#1\|}
\newcommand{\scp}[2]{\langle #1,#2\rangle}
\renewcommand{\SS}{\mathcal{S}}
\renewcommand{\AA}{\mathcal{A}}
\newcommand{\FF}{\mathcal{F}}
\newcommand{\GG}{\mathcal{G}}
\newcommand{\II}{\mathcal{I}}
\newcommand{\HH}{\mathcal{H}}
\newcommand{\KK}{\mathcal{M}}
\DeclareMathOperator{\sgn}{sgn}

\newtheorem{theorem}{Theorem}[section]
\newtheorem{lemma}[theorem]{Lemma}
\newtheorem{corollary}[theorem]{Corollary}
\newtheorem{proposition}[theorem]{Proposition}

\newtheorem{remark}[theorem]{Remark}
\newtheorem{definition}[theorem]{Definition}

\begin{document}
\title[SSN Method for Sparsity Constraints]{A Semismooth Newton Method for Tikhonov Functionals with Sparsity Constraints}

\author{R Griesse$^1$ and D A Lorenz$^2$}
\address{$^1$Johann Radon Institute for Computational and Applied Mathematics (RICAM), Austrian Academy of Sciences, Altenbergerstra{\ss}e~69, A--4040 Linz, Austria}

\address{$^2$Zentrum f\"ur Technomathematik, University of Bremen, D--28334 Bremen, Germany}

\eads{\mailto{dlorenz@math.uni-bremen.de}, \mailto{roland.griesse@oeaw.ac.at}}

\begin{abstract}
  Minimization problems in $\ell^2$ for Tikhonov functionals with sparsity constraints are considered.
  Sparsity of the solution is ensured by a weighted $\ell^1$ penalty term.
  The necessary and sufficient condition for optimality is shown to be slantly differentiable (Newton differentiable), hence a semismooth Newton method is applicable.
  Local superlinear convergence of this method is proved.
  Numerical examples are provided which show that our method compares favorably with existing approaches.
\end{abstract}






\section{Introduction}
\label{sec:Introduction}
In this work we consider the optimization problem
\begin{equation}
  \label{eq:Main_problem}
  \text{Minimize} \quad \frac{1}{2} \norm{Ku-f}^2_\HH + \sum_{k=1}^\infty w_k \, |u_k| \quad \text{over } u \in \ell^2.
\end{equation}
Here, $K:\ell^2 \to \HH$ is a linear and injective operator mapping the sequence space $\ell^2$ into a Hilbert space $\HH$, $f \in \HH$ and $w = \{w_k\}$ is a sequence satisfying $w_k \ge w_0 > 0$.

One well understood algorithm for the solution of~(\ref{eq:Main_problem}) is the so-called iterated soft-thresholding for which convergence has been proven in~\cite{daubechies2004iteratethresh}, see also \cite{combettes2005signalrecovery,bredies2007softthresholding}.
While the iterated soft-thresholding is very easy to implement it converges very slow in practice (in fact the method converges linearly but with a constant very close to one~\cite{bredies2007softthresholding}).
Another well analyzed method is the iterated hard-thresholding which converges like $\mathcal{O}(n^{-1/2})$~\cite{bredies2008harditer} (i.e.~even slower than the iterated soft-thresholding but practically it is faster in many cases).

In this article we derive an algorithm for which we prove local superlinear convergence in the infinite dimensional setting.
Our algorithm is an active set, or semismooth Newton, method and hence, the analysis is based on the notion of slant differentiability~\cite{chen2000semismooth,HintermuellerItoKunisch02}.
The semismooth Newton method is easily implementable as an active set method.
Numerical experiments show that the method is robust with respect to the choice of the initial value and that it compares favorably with existing approaches in terms of computation time.

The background for problems of type~(\ref{eq:Main_problem}) is, for example, the attempt
to solve the linear operator equation $Ku = f$ in an
infinite-dimensional Hilbert space which models the connection between
some quantity of interest $u$ and some measurements
$f$. Often, the measurements $f$ contain noise which makes the
direct inversion ill-posed and practically impossible. Thus, instead
of considering the linear equation, a regularized problem is posed for
which the solution is stable with respect to noise.
A common approach is to regularize by minimizing a
Tikhonov functional
\cite{engl1996inverseproblems,daubechies2004iteratethresh,
  resmerita2005regbanspaces}.
 A special class of these regularizations has been
of recent interest, namely of the type~(\ref{eq:Main_problem}).
These problems model the fact that the quantity of interest $u$ is
composed of a few elements, i.e.~it is sparse in some given, countable
basis.
To make this precise, let $A:\HH_1\to\HH_2$ be a bounded operator between two Hilbert spaces and let $\{\psi_k\}$ be an orthonormal basis of $\HH_1$.
Denote by $B:\ell^2\to\HH_1$ the synthesis operator
$B(u_k) = \sum_{k}u_k\psi_k$. Then the problem
\[
\min_{u\in\HH_1} \frac{1}{2} \norm{Au - f}^2_{\HH_2} +
  \sum_{k=1}^\infty w_k |\scp{u}{\psi_k}|
\]
can be rephrased as
\[
\min_{u\in \ell^2}\frac{1}{2} \norm{ABu - f}^2_{\HH_2} +
  \sum_{k=1}^\infty w_k |u_k|. 
\]
The sequence $w_k$ plays the role of the regularization parameter where
each coefficient is regularized individually.
However, for an analysis of the regularizing properties one might use $\alpha \, w_k$ instead and investigate $\alpha\to 0$.
We refer to e.g.~\cite{daubechies2004iteratethresh,ramlau2006tikhproject,lorenz2008reglp} for analysis of the regularizing properties and parameter choice rules.

Recently sparsity constraints have also appeared in the context of optimal control of PDEs \cite{Stadler07:1}.

The article is organized as follows.
In Section~\ref{sec:Optimality_Conditions} we derive a semismooth formulation for the minimization problem~(\ref{eq:Main_problem}).
Section~\ref{sec:SSN_Method}  states the algorithm and local superlinear convergence is proven.
The Section~\ref{sec:Numerical_Results} presents numerical results on the regularization of the ill-posed problems of inverse integration and deblurring and shows an application to $\ell^1$ minimization in the context of compressed sensing.

\subsection*{Notation}
For $1 \le p < \infty$, $\ell^p$ denotes the space of $p$-summable sequences with norm $\norm{u}_p = \Big( \sum_{k=1}^\infty |u_k|^p \Big)^{1/p}$, whereas $\ell^\infty$ denotes the space of bounded sequences with norm $\norm{u}_\infty = \max_{k \in \N} |u_k|$.
Recall that these spaces satisfy $\ell^p \embed \ell^q$ for $1 \le p \le q \le \infty$ and that $|u_k|\leq \norm{u}_p$ holds for any $u \in \ell^p$.
In the case $p=2$ we simply write $\norm{u}$, and $\scp{\cdot}{\cdot}$ denotes the inner product in $\ell^2$.
With $B_\rho(u)$ we denote the open ball of radius $\rho$ with respect to the norm of $\ell^2$, centered at $u$.
The operator $K^*:\HH \to \ell^2$ is the Hilbert space adjoint of $K$ and $L(X,Y)$ is the space of bounded linear operators from $X$ to $Y$.


\section{Optimality Conditions}
\label{sec:Optimality_Conditions}
In this section we are going to derive the necessary and sufficient optimality condition for the problem~(\ref{eq:Main_problem}).
It is going to be the basis for the semismooth Newton algorithm.
This condition can be derived and expressed in different ways, for example by using the classical Lagrange duality, or by using subgradient calculus.

Let us first address the conditions obtained by subgradient calculus.
To this end we introduce the so-called soft-thresholding function.
\begin{definition}
  \label{def:thresh_function}
  Let $w = \{w_k\}$ with 
  $w_k\geq w_0 > 0$ and  $1\leq p < \infty$, $1\leq q \leq \infty$.
  The \emph{soft-thresholding of $u$ with the sequence $w$} is defined as the mapping $\SS_w: \ell^p \to \ell^q$ given by
  \begin{equation}\label{eq:thresh_function}
    \SS_w(u)_k = S_{w_k}(u_k) = \max\{0,|u_k|-w_k\}\sgn(u_k).
  \end{equation}
\end{definition}

\begin{remark}\label{rem:range_soft_thresh}
  Since elements of $\ell^p$ are sequences converging to zero, the range of $\SS_w$ is $\ell^0 = \{u\in\R^\N\ :\ u_k = 0 \text{ for almost every } k\}\subset \ell^q$. 
\end{remark}

With the help of the soft-thresholding operator, we can formulate the optimality condition in a compact way.

\begin{proposition}
  \label{prop:Min_characterization_shrinkage}
  If $K:\ell^2\to \HH$ is injective, the functional
  \begin{equation}
    \label{eq:functional}
    \Psi(u) = \frac{1}{2} \norm{Ku-f}^2_\HH + \sum_{k=1}^\infty w_k |u_k|
  \end{equation}
  has a unique minimizer $\bar u \in \ell^2$.
  This minimizer is characterized by
  \begin{equation}
    \label{eq:Min_characterization_shrinkage}
    \bar u = \SS_{\gamma w}(\bar u -\gamma K^*(K\bar u -f)) \quad \text{for any } \gamma>0.
  \end{equation}
\end{proposition}
\begin{proof}
  Since $K$ is injective, $\Psi$ is strictly convex and coercive and hence,hen it has a unique minimizer. This minimizer is characterized by
  \begin{equation*}
    0 \in \partial\Psi(\bar u)
  \end{equation*}
  which is equivalent to
  \begin{equation}
    \label{eq:char_minimizer_inclusion}
    -K^*(K\bar u -f) \in \partial F(\bar u)
  \end{equation}
  where $F(u) = \sum_k w_k|u_k|$.
  Multiplying with $\gamma>0$, adding $\bar u$ to both sides and inverting $(I+\gamma\partial F)$ gives
  \begin{equation*}
    \bar u = (I+\gamma \, \partial F)^{-1}(\bar u - \gamma K^*(K \bar u-f)).
  \end{equation*}
  (Note that $(I+\gamma \, \partial F)^{-1}$ exists and is single-valued since the subgradient $\partial F$ is maximal monotone if $F$ is convex and lower semicontinuous~\cite[Proposition~32.17, Corollary~32.30]{zeidler90:_nonlin_funct_analy_applic_ii_b}.)
  A straightforward calculation shows that
  \begin{equation*}
    (I+\gamma \, \partial F)^{-1} = \SS_{\gamma w}.
  \end{equation*}
\end{proof}

From the characterization~(\ref{eq:Min_characterization_shrinkage}) and Remark~\ref{rem:range_soft_thresh} we can derive the following corollary.
\begin{corollary}
  The minimizer $\bar u$  of~(\ref{eq:functional}) is a finitely supported sequence.
\end{corollary}

For convenience we also derive the optimality condition by Lagrange duality.
We split the functional $\Psi$ according to
\[
\Psi(u) = G(Ku) + F(u)
\]
with $G(h) = \norm{h-f}^2_\HH/2$ and $F(u) = \sum_k w_k|u_k|$.
To state the dual problem, we use the dual variable $p$ which shall not be confused with exponents for $\ell^p$ spaces. 
The dual variable appears only in this section.
The dual problem of (\ref{eq:Main_problem}) is defined as
\begin{equation}
  \label{eq:Main_problem_dual}
  \text{Maximize} \quad -F^*(K^*p)-G^*(-p) \quad \text{over } p \in \HH.
\end{equation}
This can be expressed as (see the appendix)
\begin{equation*}
  \begin{aligned}
    \text{Maximize} \quad & - \frac{1}{2} \norm{p}^2_\HH + \langle p, f \rangle \quad \text{over } p \in \HH \\
    \text{subject to} \quad & |K^* p|_k \le w_k \quad \text{for all } k.
  \end{aligned}
\end{equation*}

The extremality conditions \cite[Ch.~III.4]{ekelandtemam1976convex} are:
\begin{subequations}
  \begin{align}
    \label{eq:extremality_a}
    F(u) + F^*(K^* p) - \langle K^* p, u \rangle = 0 \\
    \label{eq:extremality_b}
    G(K u) + G^*(-p) + \langle p, K u \rangle = 0.
  \end{align}
\end{subequations}
The first condition (\ref{eq:extremality_a}) yields
\begin{align*}
 & \sum_{k=1}^\infty \underbrace{w_k \, |u_k|}_{\ge 0} - \langle K^* p, u \rangle = 0 \quad \text{and} \quad |(K^*p)_k| \le w_k \\
  \Leftrightarrow \quad & u_k = 0 \quad \text{or} \quad (K^* p)_k = w_k \text{ sign } u_k = 
  \begin{cases}
    w_k, & \text{if } u_k > 0 \\
    - w_k, & \text{if } u_k < 0
  \end{cases} \\
  & \quad \text{and} \quad |(K^*p)_k| \le w_k.
\end{align*}
This condition can be written as the complementarity system
\begin{equation}
  \label{eq:extremality_a_as_complementarity_system}
  \begin{aligned}
    K^* p - w & \le 0, & \quad u^+ & \ge 0, & \quad [K^* p - w] \, u^+ & = 0 \\
   - K^* p - w & \le 0, & \quad u^- & \ge 0, & \quad [K^* p + w] \, u^- & = 0,
  \end{aligned}
\end{equation}
in a coordinatewise sense, which is in turn equivalent to
\begin{equation}
  \label{eq:extremality_a_as_semismooth_equation}
  u = \max \{ 0, u + \gamma \, (K^* p - w) \} + \min \{ 0, u + \gamma \, (K^* p + w) \}
\end{equation}
for any $\gamma > 0$.

The second condition (\ref{eq:extremality_b}) yields
\begin{equation*}
    \frac{1}{2} \norm{Ku-f}^2_\HH + \frac{1}{2} \norm{-p}^2_\HH + \langle -p, f \rangle + \langle p, Ku \rangle = 
    \frac{1}{2} \norm{Ku-f+p}^2_\HH = 0
\end{equation*}
and thus
\begin{equation}
  \label{eq:extramlity_b_evaluated}
  Ku-f+p = 0.  
\end{equation}

By plugging (\ref{eq:extramlity_b_evaluated}) into~(\ref{eq:extremality_a_as_semismooth_equation}) we end up with
\begin{equation}
  \label{eq:extremality_semismooth_system}
  u - \max \{ 0, u - \gamma \, (K^* (Ku - f) + w) \} - \min \{ 0, u - \gamma \, (K^* (Ku - f) - w) \} = 0,
\end{equation}
which is just another way to express (\ref{eq:Min_characterization_shrinkage}).

\begin{remark}
The usual characterization $0\in\partial\Psi(\bar u)$ of the unique minimizer $\bar u$ of~(\ref{eq:Main_problem}) is diffucult to handle for numerical algorithms because it is a nonsmooth inclusion.
One attempt to tackle the problem is by interior point regularization as proposed in~\cite{kim2007l1ls}.
This, however, introduces additional nonlinearities into the problem.
By contrast, our algorithm is based on the necessary and sufficient condition~(\ref{eq:extremality_semismooth_system}).
As we shall prove in the following section,~(\ref{eq:extremality_semismooth_system}) is a semismooth equation in $\ell^2$, so that Newton's method can be applied.
\end{remark}


\section{Semismooth Newton Method}
\label{sec:SSN_Method}
The previous section has shown that we can solve the minimization problem~(\ref{eq:Main_problem}) by solving the equation (\ref{eq:Min_characterization_shrinkage}) or (\ref{eq:extremality_semismooth_system}), or briefly
\begin{equation}
  \label{eq:opti_condition_F}
  \FF(u) = u - \SS_{\gamma w}(u - \gamma K^*(Ku-f)) = 0,
\end{equation}
for some $\gamma>0$.

This is an operator equation in the space $\ell^2$, involving the non-differentiable $\max$ and $\min$ operations.
Optimality conditions of this form frequently also occur in the context of optimal control problems for partial differential equations, in the presence of control constraints.
Then (\ref{eq:opti_condition_F}) is considered in $L^p$ function spaces, and it is known that the $\max$ operation, i.e., $u\mapsto \max\{0,u\}$, is so-called Newton or slantly differentiable from $L^p$ to $L^q$ for $1 \le q \le p \le \infty$, see \cite[Theorem~2.6]{chen2000semismooth} in view of its Lipschitz continuity.
In the presence of a norm gap $1 \le q < p \le \infty$, the generalized derivative, or slanting function, can be chosen as an indicator function, see \cite[Proposition~4.1]{HintermuellerItoKunisch02}.
This allows for the interpretation of the generalized Newton method as a so-called active set method.
This norm gap is made up for in the context of partial differential equation because $K$ and $K^*$ are solution operators which provide the necessary smoothing.

It turns out that the behavior of the $\max$ and $\min$ operations is more intricate than in function space.
Again, it follows from the Lipschitz continuity of $u \mapsto \max\{0,u\}$ from $\ell^p$ to $\ell^q$ that slant differentiability holds \cite[Theorem~2.6]{chen2000semismooth} for $1 \le p, q \le \infty$.
However, we are not aware of any simple slanting function even with norm gap which can be algoithmically exploited, see Remark~\ref{rem:no_slant_func}.
It may be surprising that nonetheless, the soft-thresholding operator $\SS_w$ and thus equation (\ref{eq:opti_condition_F}) are slantly differentiable and admit a simple slanting function between \emph{any} pair of $\ell^p$, $\ell^q$ spaces, see Proposition~\ref{prop:slantderivative_shrinkage}.
This allows us to apply a generalized Newton's method to solve (\ref{eq:opti_condition_F}), which takes the form of an active set method.

\subsection{Semismoothness of the optimality condition}
\label{sec:semism-optim-cond}

The concept of slant or Newton differentiability is closely related to the notion of semismoothness \cite{chen2000semismooth,HintermuellerItoKunisch02,Ulbrich03}, and we will use the terms interchangeably.

\begin{definition}
  Let $X$ and $Y$ be Banach spaces and $D\subset X$ be an open subset.
  A mapping $\FF:D\to Y$ is called \emph{Newton (or slantly) differentiable in $x\in D$} if there exists a family of mappings $\GG:D\to L(X,Y)$ such that
  \begin{equation}
    \label{eq:def_slant_diff}
    \lim_{h\to 0} \frac{\norm{\FF(x+h) - \FF(x) - \GG(x+h) \, h}_Y}{\norm{h}_X}= 0.
  \end{equation}
  The function $\GG$ is called a \emph{generalized derivative (or slanting function) for $\FF$ in $x$}. 
\end{definition}

It is shown in~\cite{chen2000semismooth} that any Lipschitz continuous function is Newton differentiable.
However, this is only of little help algorithmically unless there is a generalized deriative $\GG(u)$ of~(\ref{eq:opti_condition_F}) which is easily invertible. 

\begin{remark}
  \label{rem:no_slant_func}
  A natural candidate for a generalized derivative $\GG$ of the function $\FF(u) = \max(0,u)$ is
  \[\GG(u)(h)_k  = 
  \begin{cases}
    h_k & , u_k>0\\
    \delta h_k &, u_k=0\\
    0 & , u_k<0 
  \end{cases}
  \quad\text{for any } \delta\in\R.
  \]
  We are going to show that this $\GG$ can not serve as a generalized derivative of $\FF:\ell^p\to \ell^q$ for any $p\in[1,\infty[$ and $1\leq q\leq \infty$.
  We consider a point $u\in \ell^p$ for which the set $\{n\ |\ u_n\not=0\}$ is infinite and 
  take a special sequence of $h^n\in \ell^p$, namely
  \begin{equation*}
    h_k^n = 
    \begin{cases}
      0 & \text{for } k \neq n \\
      -2 u_k & \text{for } k = n.
    \end{cases}
  \end{equation*}
  Hence, we have $\norm{h^n}_p = 2|u_n|\to 0$ for $n\to\infty$.
  It is an easy calculation to see that
  \begin{equation*}
    \frac{\norm{\max\{u+h^n,0\} - \max\{u,0\} - G(u+h^n)h^n}_q}{\norm{h^n}_p} = \frac12\ \text{ for all } n \text{ with } u_n\not=0.
  \end{equation*}
\end{remark}

The following proposition shows that the thresholding operator~(\ref{eq:thresh_function}) is Newton differentiable and that a function similar to $\GG$ serves as a generalized derivative.

\begin{proposition}
  \label{prop:slantderivative_shrinkage}
  The mapping $\SS_w:\ell^p \to \ell^q$  from Definition~\ref{def:thresh_function} is Newton differentiable for any $1\leq p < \infty$, $1\leq q \leq \infty$.
  A generalized derivative is given by
  \begin{equation*}
    (\GG(u) \, v)_k = 
    \begin{cases}
      v_k & \text{for } |u_k|>w_k\\
      0   & \text{for } |u_k|\leq w_k.
    \end{cases}
  \end{equation*}
\end{proposition}
\begin{proof}
  Without loss of generality we may assume $\norm{h}_p < \tfrac{w_0}{2}$
  and hence $|h_k| < \tfrac{w_0}{2}$.
  Since $u\in \ell^p$ with $p< \infty$ there exists $k_0$ such that $|u_k| < \tfrac{w_0}{2}$ for $k>k_0$.
  We estimate
  \begin{multline*}
    \norm{\SS_w(u+h) - \SS_w(u) - \GG(u+h)(h)}^q_q \\
    \begin{aligned}
      & = \sum_{k=1}^\infty \big|S_{w_k}(u_k+h_k) - S_{w_k}(u_k) - \GG(u+h)(h)_k \big|^q \\
      & = \sum_{{k \leq k_0 \atop |u_k|\not= w_k}} \big|S_{w_k}(u_k+h_k) - S_{w_k}(u_k) - \GG(u+h)(h)_k \big|^q.
    \end{aligned}
  \end{multline*}
  It is easy to check that the above sum is zero for
  \begin{equation*}
    \norm{h}_p < \min\{ \bigl| |u_k| - w_k\bigr|\ :\ k\leq k_0 \text{ and } |u_k|\not= w_k\}
  \end{equation*}
  because $|h_k|\leq \norm{h}_p$ holds.
  It follows that
  \begin{equation*}
    \frac{\norm{\SS_w(u+h) - \SS_w(u) - \GG(u+h)(h)}_q}{\norm{h}_p} = 0
  \end{equation*}
  for $\norm{h}_p$ small enough, which proves Newton differentiability.
\end{proof}
\begin{remark}
  In matrix notation we can express the generalized derivative $\GG(u)$ as
  \begin{equation*}
    \GG(u) = 
    \begin{pmatrix}
      I_\AA & 0 \\
      0     & 0
    \end{pmatrix}
  \end{equation*}
  where $\AA = \{k\in\N\ :\ |u_k|>w_k\}$.
\end{remark}

To calculate a generalized derivative for the mapping $\FF$ in~(\ref{eq:opti_condition_F}), we prove a chain rule for the generalized derivative.

\begin{lemma}
  \label{lem:slantderivative_chainrule}
  Let $S:X\to Y$ be Newton differentiable, $A\in L(X,X)$ and $y\in X$.
  Let furthermore $\GG$ be a generalized derivative of $S$.
  Define $T(u) = S(Au + y)$.
  Then $H(u) = \GG(Au+y)A$ is a generalized derivative of $T$.
\end{lemma}
\begin{proof}
  It holds
  \begin{multline*}
    \frac{\norm{T(u+h) - T(u) - H(u+h)\, h}}{\norm{h}}\\
     = \frac{\norm{S(Au + Ah +y) - S(Au+y) - \GG(Au + Ah + y)Ah}}{\norm{Ah}}\frac{\norm{Ah}}{\norm{h}}.
  \end{multline*}
  The right hand side converges to zero because $\GG$ is a generalized derivative of $S$ in $Au + y$ in the direction $Ah$, and $A$ is bounded.
\end{proof}

In order to specify a generalized derivative of $\FF$, we introduce the active and the inactive sets.
For the sake of simplicity we will restrict ourself to the case $\FF:\ell^2\to \ell^2$ in the following
\begin{definition}
  \label{def:act_inact_sets}
  For $u\in \ell^2$, the active set $\AA(u)$ and the inactive set $\II(u)$ are given by
  \begin{equation*}
    \begin{aligned}
      \AA(u) & = \{k\in\N : |u - \gamma K^*(Ku-f)|_k > \gamma \, w_k\} \\
      \II(u) & = \{k\in\N : |u - \gamma K^*(Ku-f)|_k \leq \gamma \, w_k\}.
    \end{aligned}
  \end{equation*}
  Whenever the active and inactive sets correspond to an iterate $u^n$, we will denote them by $\AA_n$ and $\II_n$, respectively.
  We will drop the subscript or the argument if no ambiguity can occur.
\end{definition}

We are now in the position to calculate a generalized derivative of $\FF$.
\begin{proposition}
  \label{prop:slantderivative}
  The mapping $\FF: \ell^2 \to \ell^2$,
  \begin{equation*}
    \label{eq:shrinkage}
    \FF(u) =  u - \SS_{\gamma w}(u - \gamma K^*(Ku - f))
  \end{equation*}
  is Newton differentiable.
  Denote the active and inactive set $\AA$ and $\II$ as in Definition~\ref{def:act_inact_sets}
  and split the operator $K^* K$ according to
  \begin{equation*}
    \label{eq:opsplit}
    K^* K = 
    \begin{pmatrix}
      \KK_{\AA\AA} & \KK_{\AA\II}\\
      \KK_{\II\AA} & \KK_{\II\II}
    \end{pmatrix}.
  \end{equation*}
  Then a generalized derivative is given by
  \begin{equation}\label{eq:slant_func}
    \GG(u) = 
      \begin{pmatrix}
        0 & 0\\
        0 & I_\II
      \end{pmatrix}
       + 
      \begin{pmatrix}
       I_\AA & 0\\
       0     & 0
      \end{pmatrix}
      (\gamma K^*K)
      = 
      \begin{pmatrix}
        \gamma \KK_{\AA\AA} & \gamma \KK_{\AA\II}\\
        0                 & I_{\II}
      \end{pmatrix}.
  \end{equation}
\end{proposition}
\begin{proof}
  The claim follows from the sum rule for the generaized derivative and from Proposition~\ref{prop:slantderivative_shrinkage} and Lemma~\ref{lem:slantderivative_chainrule} with $S= \SS_{\gamma w}$, $A = I-\gamma K^*K$, $y = \gamma K^*f$.
\end{proof}
\begin{remark}
  \label{rem:Finiteness_of_A}
  Note that for any $u \in \ell^2$, the active set $\AA$ is always finite, since $u - \gamma K^*(Ku-f) \in \ell^2$ holds and thus $|u - \gamma K^*(Ku-f)|_k \to 0$ for $k \to \infty.$
\end{remark}

\subsection{Semismooth Newton method}

The semismooth or generalized Newton method for the solution of (\ref{eq:opti_condition_F}) can be stated as the iteration
\begin{equation}
  \label{eq:ssn_method_general}
  u^{n+1} = u^n - \GG(u^n)^{-1}\FF(u^n),
\end{equation}
where $\GG$ is a generalized derivative of $\FF$.
We use the generalized derivative $\GG$ given by~(\ref{eq:slant_func}).
Naturally, the semismooth Newton method can be interpreted as an active set method, and we state it as Algorithm~\ref{alg:SSN}.

\begin{algorithm}
  \caption{Semismooth Newton method for the solution of (\ref{eq:opti_condition_F}).}
  \label{alg:SSN}
  \begin{algorithmic}[1]
    \STATE Initialize $u^0$, choose $\gamma > 0$, set $n := 0$ and done := false
    \WHILE{$n < n_{max}$ and not done}
    \STATE Calculate the active and inactive sets:
    \begin{align*}
      \AA & = \{k\in\N : |u^n - \gamma K^*(Ku^n-f)|_k > \gamma w_k\} \\
      \II & = \{k\in\N : |u^n - \gamma K^*(Ku^n-f)|_k \leq \gamma w_k\}.
    \end{align*}
    \STATE Compute the residual
    \begin{align*}
      r^n = \FF(u^n) = u^n - \SS_{\gamma w}(u^n - \gamma K^*(Ku^n-f)).
    \end{align*}
    \IF{$\norm{r^n} \le \varepsilon$}
    \STATE done := true
    \ELSE
    \STATE Calculate the Newton update by solving
    \begin{equation*}
      \begin{pmatrix}
        \gamma \KK_{\AA\AA} & \gamma \KK_{\AA\II}\\
        0                 & I_{\II}
      \end{pmatrix}
      \begin{pmatrix}
        \delta u_{\AA}\\
        \delta u_{\II}
      \end{pmatrix}
      = - 
      \begin{pmatrix}
        r^n_\AA\\
        r^n_\II
      \end{pmatrix}
    \end{equation*}
    \STATE Update $u^{n+1} := u^n + \delta u$
    \STATE Set $n := n+1$
    \ENDIF
    \ENDWHILE
  \end{algorithmic}
\end{algorithm}
 
\begin{remark}
  \begin{enumerate}
  \item Algorithm~\ref{alg:SSN} is the generalized Newton method (\ref{eq:ssn_method_general}).
    The unique solvability in step~8 is shown in Proposition~\ref{prop:Boundedness_of_the_inverse} below.

  \item Given an initial iterate $u^0 \in \ell^2$, the algorithm is well-defined, and all iterates remain in $\ell^2$.
    We refer again to Proposition~\ref{prop:Boundedness_of_the_inverse} below.

  \item   At the end of step~10, the iterate $u^{n+1}$ satisfies $u^{n+1}_{\II_n} = 0$.
    Note that $r^n_{\II_n} = u^n_{\II_n}$ holds which implies $\delta u_{\II_n} = - u^n_{\II_n}$.

  \end{enumerate}
\end{remark}
Note that $(iii)$ implies that all iterates $u^n$ ($n \ge 1$) of Algorithm~\ref{alg:SSN} are finitely supported sequences.
However, $K^*(K u^n - f)$ is in general not finitely supported, and hence in a practical implentation, this term will be truncated after a number of entries.

\subsection{Active set method}
\label{sec:active-set-method}
One may set up Algorithm~\ref{alg:SSN} equivalently as an active set method.
This can be seen by a closer analysis of the Newton step (step~8 and~9 in Algorithm~\ref{alg:SSN}):
\begin{align*}
  u^{n+1} & = u^n -
  \begin{pmatrix}
    \tfrac{1}{\gamma} \KK_{\AA\AA}^{-1} & -\KK_{\AA\AA}^{-1} \, \KK_{\AA\II}\\
    0                & I_\II
  \end{pmatrix}
  \begin{pmatrix}
    u^n - \SS_{\gamma w}(u^n - \gamma K^*(Ku^n-f))
  \end{pmatrix}\\
  & = u^n - 
  \begin{pmatrix}
    \tfrac{1}{\gamma} \KK_{\AA\AA}^{-1} & -\KK_{\AA\AA}^{-1} \, \KK_{\AA\II}\\
    0                & I_\II
  \end{pmatrix}
  \begin{pmatrix}
    \gamma \, [K^*(Ku^n-f)]_\AA \pm w_\AA)\\
    u^n_\II
  \end{pmatrix}\\
  & = 
  \begin{pmatrix}
    u^n_\AA - \KK_{\AA\AA}^{-1}\bigl( [K^*(Ku^n-f)]_\AA \pm w_\AA -\KK_{\AA\II} \, u^n_\II\bigr)\\
    0               
  \end{pmatrix}\\
  & = 
  \begin{pmatrix}
    \KK_{\AA\AA}^{-1}(K^*f \pm w)|_\AA\\
    0               
  \end{pmatrix}
\end{align*}
The sign of $w$ depends of the sign of $u^n - \gamma K^* (Ku^n-f)$.
Hence, instead of calculating the Newton update in step~8, one may set $u^{n+1}_{\II} := 0$ and solve $\KK_{\AA\AA} u^{n+1}_{\AA} = (K^* f \pm w)_{|\AA}$.
This shows that the subsequent iterate $u^{n+1}$ depends on the current iterate $u^n$ solely through the active set $\AA$.
As a consequence, differing values of $u^n$ can lead to the same next iterate $u^{n+1}$.

For completeness, we state the active set method as Algorithm~\ref{alg:active_set}.
Note that the algorithm is initialized with an active set $\AA$ instead of $u^0$.
\begin{algorithm}
  \caption{Active set method for the solution of (\ref{eq:opti_condition_F}).}
  \label{alg:active_set}
  \begin{algorithmic}[1]
    \STATE Initialize $\AA_0^+$, $\AA_0^-$, choose $\gamma > 0$, set $n := 0$ and done := false
    \STATE Set $\AA_0 = \AA_0^+\cup \AA_0^-$, $\II_0 = \N\setminus \AA_0$
    \STATE Set the signs of the weights:
    \begin{equation*}
      s^0_k =
      \begin{cases}
        1, & k \in \AA_0^+\\
        0, & k \in \II_0 \\
        -1, & k \in \AA_0^-
      \end{cases}
    \end{equation*}
    \WHILE{$n < n_{max}$ and not done}
    \STATE Set $u^n_{\II_n} = 0 $ and calculate $u^n_{\AA_n}$ by solving
    \begin{equation*}
     \KK_{\AA_n\AA_n} u^{n}_{\AA_n} = (K^* f + s^n w)_{|\AA_n}
    \end{equation*}
    \STATE Calculate the new active sets:
    \begin{align*}
      \AA_{n+1}^+ & = \{k\in\N : [u^{n} - \gamma K^*(Ku^{n}-f)]_k > \gamma w_k\} \\
      \AA_{n+1}^- & = \{k\in\N : [u^{n} - \gamma K^*(Ku^{n}-f)]_k < -\gamma w_k\} \\
      \II_{n+1} & = \{k\in\N : |u^{n} - \gamma K^*(Ku^{n}-f)|_k \leq \gamma w_k\}.
    \end{align*}
    \STATE Set the signs of the weights:
    \begin{equation*}
      s^{n+1}_k =
      \begin{cases}
        1, & k \in \AA_{n+1}^+\\
        0, & k \in \II_{n+1} \\
        -1, & k \in \AA_{n+1}^-
      \end{cases}
    \end{equation*}
    \IF {$s^{n+1} = s^{n}$}
      \STATE done := true
    \ENDIF
    \STATE Set $n:= n+1$
    \ENDWHILE
  \end{algorithmic}
\end{algorithm}

In this setting, the stopping criterion is coincidence of the active sets in consecutive iterations---other choices are also possible.
In the numerical examples in Section~\ref{sec:Numerical_Results} we chose the norm of the residual because a sudden drop of the residual norm occured before the minimizer was identified.

\subsection{Local convergence of the semismooth Newton method}
\label{sec:conv-semism-newt}

The local superlinear convergence of the semismooth Newton method (Algorithm~\ref{alg:SSN}) hinges upon the uniform boundedness of $\GG(u^n)^{-1}$ during the iteration.
\begin{proposition}
  \label{prop:Iterates_are_uniformly_finitely_supported}
  There exists $k_0 \in \N$ and $\rho > 0$ such that $\norm{u-\overline{u}}_2 < \rho$ implies that 
  \begin{equation*}
    \AA(u) \subset [1,k_0].
  \end{equation*}
  Moreover, $k_0$ and $\rho$ depend only on $\gamma$, $\overline{u}$, $\norm{K^*K}$, $\norm{K^*f}$, and $w_0$.
\end{proposition}
\begin{proof}
  The triangle inequality implies
  \begin{equation}
    \label{eq:Iterates_are_uniformly_finitely_supported_1}
    |u - \gamma K^*(Ku - f)|_k
    \le |\overline{u} - \gamma K^*(K\overline{u} - f)|_k
    + |u - \overline{u} - \gamma K^*K(u - \overline{u})|_k.
  \end{equation}
  The first term can be estimated by
  \begin{equation*}
    |\overline{u} - \gamma K^*(K\overline{u} - f)|_k
    \le |\overline{u} |_k + \gamma \, |K^*K \overline{u}|_k + \gamma \, |K^* f|_k.
  \end{equation*}
  Since $\overline{u}$, $K^*K \overline{u}$ and $K^* f$ are in $\ell^2$, the right hand side converges to 0 as $k \to \infty$.
  In particular, there exists $k_0$, depending only on the named quantities, such that 
  \begin{equation}
    \label{eq:Iterates_are_uniformly_finitely_supported_2}
    |\overline{u} - \gamma K^*(K\overline{u} - f)|_k
    \le \gamma \, w_0 / 2 \quad \text{for all } k \ge k_0.
  \end{equation}
  The second term in (\ref{eq:Iterates_are_uniformly_finitely_supported_1}) can be estimated by
  \begin{multline*}
    |u - \overline{u} - \gamma K^*K(u - \overline{u})|_k 
    \le |u - \overline{u}|_k + \gamma \, |K^*K(u - \overline{u})|_k \\
    \le \norm{u - \overline{u}} + \gamma \, \norm{K^*K(u - \overline{u})} 
    \le (1 + \gamma \, \norm{K^*K}) \, \norm{u - \overline{u}}.
  \end{multline*}
  Hence there exists $\rho >0$, depending only on the named quantities, such that 
  \begin{equation}
    \label{eq:Iterates_are_uniformly_finitely_supported_3}
    |u - \overline{u} - \gamma K^*K(u - \overline{u})|_k \le \gamma \, w_0 / 2 \quad \text{for all } k \in \N.
  \end{equation}
  Combining (\ref{eq:Iterates_are_uniformly_finitely_supported_1})--(\ref{eq:Iterates_are_uniformly_finitely_supported_3}) proves the claim.
\end{proof}
At this point, we cannot yet conclude that the active sets remain uniformly bounded during the iteration of Algorithm~\ref{alg:SSN}, since it is not evident whether the iterates will remain in a suitable $\rho$-neighborhood of $\overline{u}$.
\begin{proposition}
  \label{prop:Boundedness_of_the_inverse}
  The generalized derivative $\GG$, given by (\ref{eq:slant_func}), is boundedly invertible from $\ell^2$ into $\ell^2$.
  Moreover, the norm of $\GG(u)^{-1}$ can be estimated by
  \begin{equation*}
    \norm{\GG(u)^{-1}}
    \le \norm{\KK_{\AA\AA}^{-1}} \big( \tfrac{1}{\gamma} + \norm{\KK_{\AA\II}} \big) + 1,
  \end{equation*}
  where $\AA$ and $\II$ are the active and inactive sets at $u$, see Definition~\ref{def:act_inact_sets}.
\end{proposition}
\begin{proof}
  Let $u,r \in \ell^2$ and consider the equation $\GG(u) \, \delta u = r$, i.e., 
  \begin{equation*}
    \begin{pmatrix}
      \gamma \KK_{\AA\AA} & \gamma \KK_{\AA\II} \\
      0                 & I_{\II}
    \end{pmatrix}
    \begin{pmatrix}
      \delta u_\AA \\
      \delta u_\II
    \end{pmatrix}
    = 
    \begin{pmatrix}
      r_\AA \\ r_\II
    \end{pmatrix}.
  \end{equation*}
  Necessarily, $\delta u_\II = r_\II$ holds, which implies $\delta u_\II \in \ell^2$.
  It remains to solve
  \begin{equation}
    \label{eq:Inverse_on_A}
    \gamma \KK_{\AA\AA} \, \delta u_\AA = r_\AA - \gamma \KK_{\AA\II} \, \delta u_\II.
  \end{equation}
  The right hand side is an element of $\ell^2$.
  Moreover, $\KK_{\AA\AA}$ is injective.
  (We rewrite $\KK_{\AA\AA} = P_\AA K^* K P_\AA = (K P_\AA)^*K P_\AA$, where  $P_\AA$ the projection of $\ell^2$ onto the active set.
  Then $\KK_{\AA\AA}u=0$ implies $\norm{K P_\AA u}^2 = \scp{u}{\KK_{\AA\AA}u} = 0$, and hence $u_\AA = 0$ since $K$ is injective.)
  By Remark~\ref{rem:Finiteness_of_A}, the active set is finite, and thus $\KK_{\AA\AA}$ is an injective operator on a finite dimensional space, hence it is also surjective.
  We conclude that (\ref{eq:Inverse_on_A}) has a unique solution $\delta u_\AA \in \ell^2$, hence $\GG(u)^{-1}: \ell^2 \to \ell^2$ exists.
  
  The norm estimate follows from 
  \begin{equation*}
    \begin{aligned}
      \norm{\GG(u)^{-1} r} 
      & = \left\| \begin{pmatrix} \tfrac{1}{\gamma}\KK_{\AA\AA}^{-1} & -\KK_{\AA\AA}^{-1} \KK_{\AA\II} \\ 0 & I_\II \end{pmatrix} \begin{pmatrix} r_\AA \\ r_\II \end{pmatrix} \right\| \\
      & \le \tfrac{1}{\gamma} \norm{\KK_{\AA\AA}^{-1}} \norm{r_\AA} + \norm{\KK_{\AA\AA}^{-1}} \norm{\KK_{\AA\II}} \norm{r_\II} + \norm{r_\II} \\
      & \le \left( \norm{\KK_{\AA\AA}^{-1}} \big( \tfrac{1}{\gamma} + \norm{\KK_{\AA\II}} \big) + 1 \right) \norm{r}.
    \end{aligned}
  \end{equation*}
\end{proof}

\begin{corollary}
  \label{cor:Boundedness_of_the_inverse_2}
  Let $k_0 \in \N$ and $\rho > 0$ be as in Proposition~\ref{prop:Iterates_are_uniformly_finitely_supported}.
  Then $\GG(u)^{-1}$ is uniformly bounded on $B_\rho(\overline{u})$.
\end{corollary}
\begin{proof}
  Let $u \in \ell^2$ such that $\norm{u-\overline{u}} < \rho$.
  By Proposition~\ref{prop:Iterates_are_uniformly_finitely_supported}, the active set satisfies $\AA(u) \subset [1,k_0]$.
  Our plan is to show that $\norm{\GG(u)^{-1}}$ indeed depends only on $k_0$.  Indeed, we define
  \begin{equation*}
    C(k_0) := \max_{\emptyset \neq \AA \subset [1,k_0]} \norm{\KK_{\AA\AA}^{-1}} > 0.
  \end{equation*}
  Note that for every $\AA \subset [1,k_0]$, $\AA \neq \emptyset$, $\KK_{\AA\AA}$ is boundedly invertible, hence $C(k_0)$ is the maximum of finitely many positive numbers.
  Moreover, $\KK_{\AA\II}$ is obtained from $K^*K$ by restriction and extension, hence $\norm{\KK_{\AA\II}} \le \norm{K^*K}$ holds, for all choices of $\AA$ and $\II$.
  From Proposition~\ref{prop:Boundedness_of_the_inverse}, we conclude that
  \begin{equation*}
    \norm{\GG(u)^{-1}} \le C(k_0)\big( \tfrac{1}{\gamma} + \norm{K^*K} \big) + 1.
  \end{equation*}
\end{proof}

We may now combine the results above to argue the local superlinear convergence of Algorithm~\ref{alg:SSN}.
\begin{theorem}
  \label{theorem:Local_superlinear_convergence}
  There exists a radius $r \in (0,\rho]$ such that $\norm{u^0 - \overline{u}} < r$ implies that all iterates of Algorithm~\ref{alg:SSN} satisfy $\norm{u^n-\overline{u}} < r$, and $u^n \to \overline{u}$ superlinearly.
\end{theorem}
\begin{proof}
  By Corollary~\ref{cor:Boundedness_of_the_inverse_2}, the inverse of the generalized derivative, $\GG(u)^{-1}$, remains uniformly bounded in $B_\rho(\overline{u})$.
  The result is then a standard conclusion for generalized Newton methods, see \cite[Remark~2.7]{Chen97}, or \cite[Theorem~1.1]{HintermuellerItoKunisch02}.
\end{proof}

\begin{remark}
  \label{remark:Local_superlinear_convergence}
  \begin{enumerate}
  \item The neighborhood in which superlinear convergence occurs is unknown and may be small.
    The global convergence behavior of the algorithm thus deserves further investigation.
    The numerical experiments in the following section suggest that the choice of $\gamma$ is essential in achieving convergence from a bad initial guess. 
    For a related problem in Hilbert spaces with a standard Tikhonov regularization term $\norm{u}^2$, global convergence without rates was proved in \cite{RoeschKunisch02}.
    
  \item The proof of Proposition~\ref{prop:slantderivative_shrinkage} together with the chain rule (Lemma~\ref{lem:slantderivative_chainrule}) shows that the remainder
    \begin{equation*}
      \FF(u^n) - \FF(\overline{u}) - \GG(\overline{u}) (u^n - \overline{u})
    \end{equation*}
    is exactly zero for sufficiently small $\norm{u^n-\overline{u}}$.
    Hence we expect convergence in one step sufficiently close to the solution, which is confirmed by the numerical results in the following section.
  \end{enumerate}
\end{remark}

\begin{remark}
  The assumption on the injectivity of $K$  may be relaxed.
  The proof of Corollary~\ref{cor:Boundedness_of_the_inverse_2} shows that we only need that all submatrices $\KK_{\AA\AA}$ for $\AA\subset [1,k_0]$ are invertible.
  Hence, local superlinear convergence can also be proved when the $K$ satisfies the \emph{finite basis injectivity} (FBI) property~\cite{bredies2007softthresholding}.
  The FBI property states that any submatrix of $K$ consisting of a finite number of columns is injective.
  The FBI property is related to the so-called \emph{restricted isometry property} (RIP), see e.g.~\cite{baraniuk2007riprandommatrices}, which plays an important role in the analysis of minimizers of $\ell^1$ constrained problems in the theory of compressed sensing~\cite{candes2005decoding}.
\end{remark}


\section{Numerical Results}
\label{sec:Numerical_Results}
In this section we present results of numerical experiments illustrating the performance of the semismooth Newton (SSN) method.
We implemented the SSN method in MATLAB\textsuperscript{\textregistered} and made experiments on a desktop PC with an AMD Athlon\textsuperscript{\texttrademark} 64 X2.
Moreover, we are going to compare the SSN method to other state-of-the-art methods for the minimization of $\ell^1$ constrained problems, namely the GPSR methods~\cite{figueiredo2007gradproj} and the \verb|l1_ls| toolbox~\cite{kim2007l1ls}
where we used the freely available MATLAB\textsuperscript{\textregistered} implementations of these methods.
The GPSR method is based on gradient projection method with Barzilai-Borwein stepsizes and is known to converge $r$-linearly.
The \verb!l1_ls! method is a truncated Newton interior point method which is applied directly to the objective functional (note that we apply a Newton method to a reformulated optimality condition).
In addition we included the widely used iterative soft-thresholding from~\cite{daubechies2004iteratethresh} in our comparison.
Note that both the GPSR and the \verb!l1_ls! methods are set up and analyzed in a finite dimensional setting while our analysis on the SSN as well as the analysis for the iterative soft-thresholding is infinite dimensional.

\subsection{Inverse integration}
\label{sec:inverse-integration}
The problem under consideration is the classical ill-posed problem of inverse integration (or differentiation~\cite{hanke2001invproblight,schoepfer2006illposedbanach,bredies2008harditer}), i.e.~the operator $K:L^2([0,1])\to L^2([0,1])$ given by
\[ Ku(t) = \int_0^t u(s) \, ds,\quad t\in[0,1].\]
The data $f$ is given as $(f(t_k))_{k=1,\dots,N}$ with $t_k =
\tfrac{1}{N}k$. We discretized the operator $K$ by the matrix
\[K = \frac{1}{N}\begin{pmatrix}
                        1 & 0 & \dots & 0\\
                        \vdots & \ddots & \ddots & \vdots\\
                        \vdots & & \ddots & 0\\
                        1 & \dots & \dots & 1 
                    \end{pmatrix},\
  K:\R^N\to\R^N.\]
The minimization problem reads
\begin{equation}
\label{eq:min_prob_inv_int}
\min_{u\in\R^N} \ \frac{1}{2} \sum_{i=1}^N \bigl((Ku)_i - f^\delta_{i} \bigr)^2 +
\sum_{k=1}^{N} w_k |u_k|.
\end{equation}
One can check easily that the SSN method is also applicable in finite dimensions and hence, this minimization problem can be treated by the SSN method.
The discussion of the SSN method in infinite dimension provides us with results
which are independent of the dimension $N$, i.e.~the algorithm scales well.

\begin{figure}[htb]
  \centering
  \includegraphics[width=.32\textwidth]{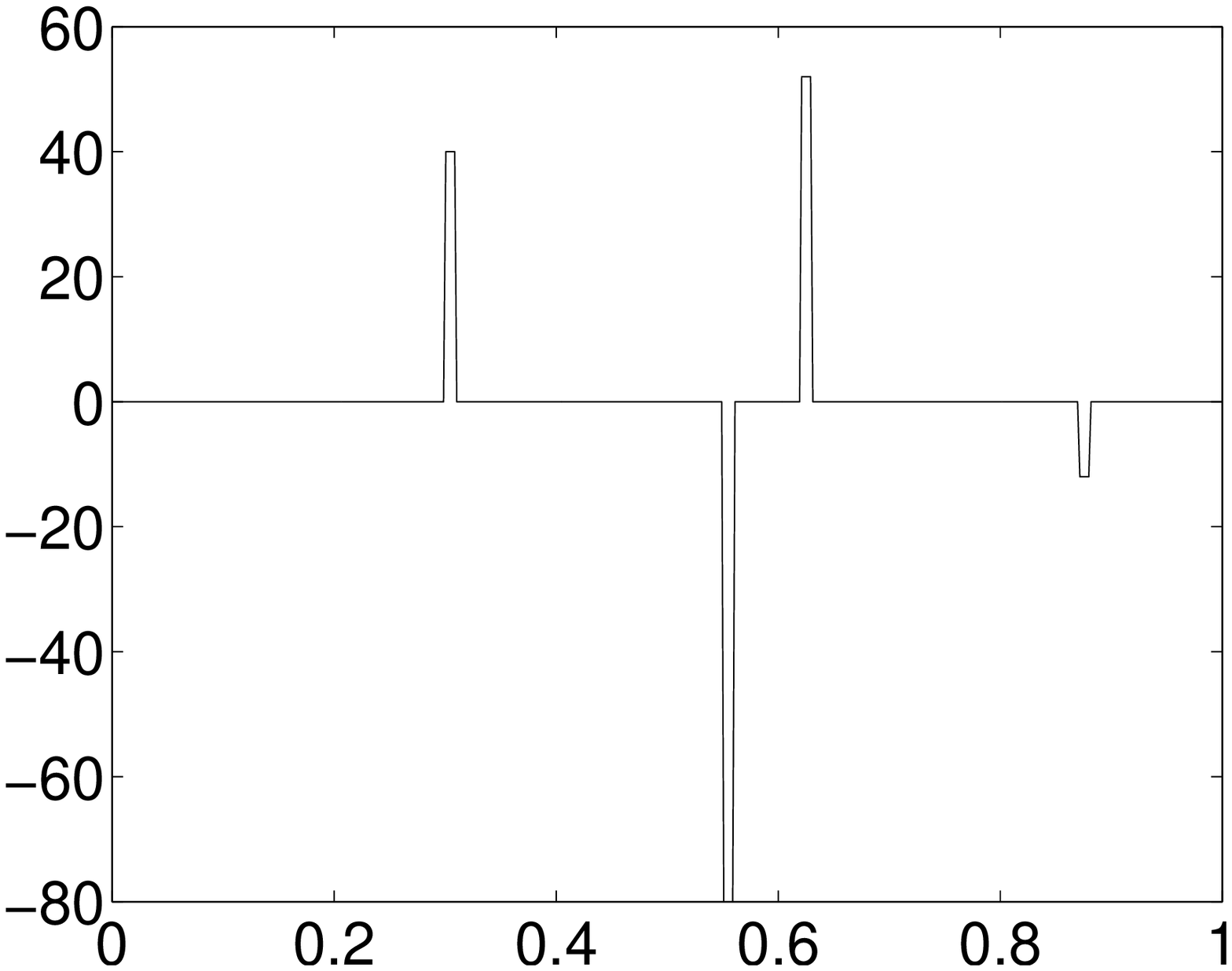}
  \includegraphics[width=.32\textwidth]{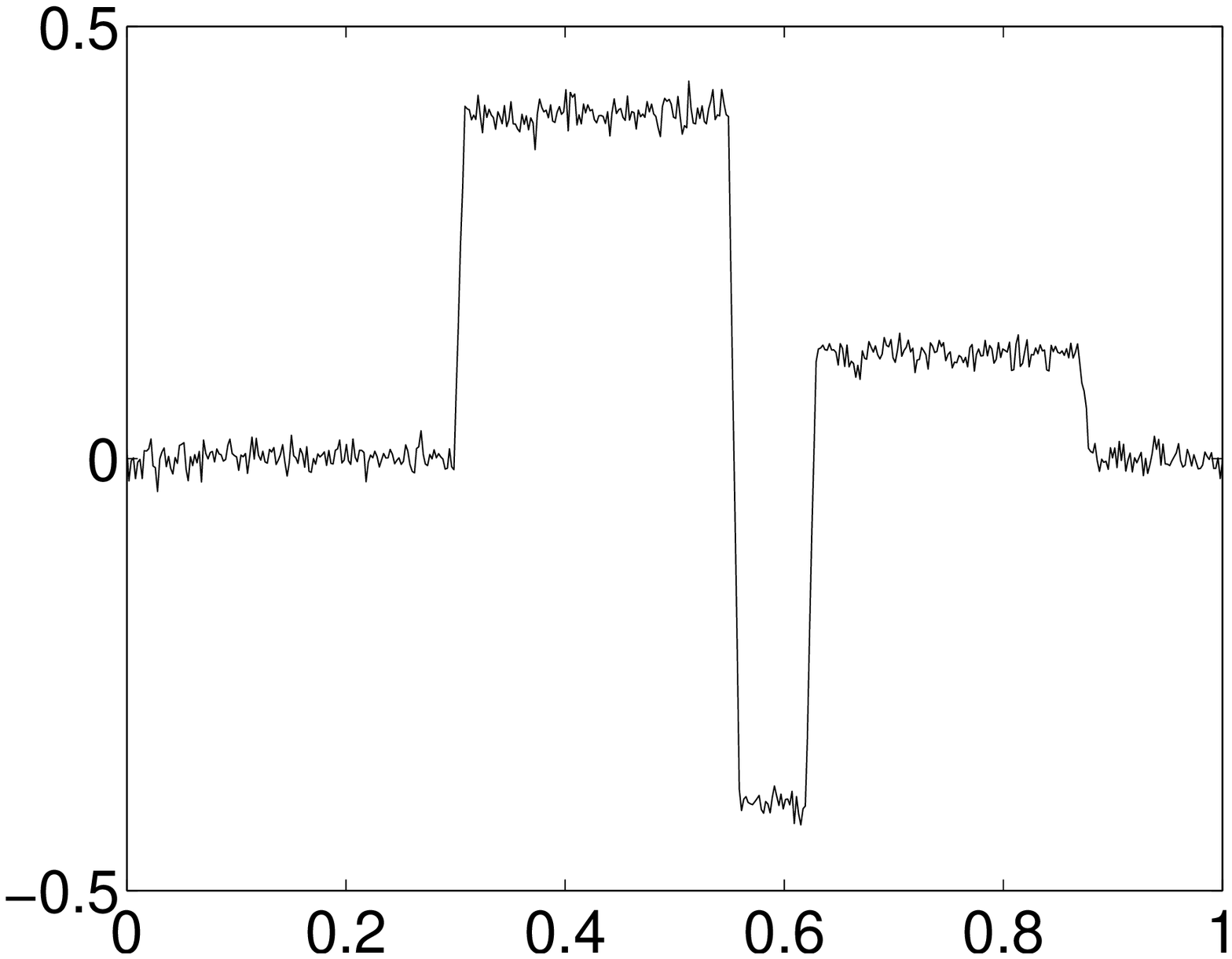}
  \includegraphics[width=.32\textwidth]{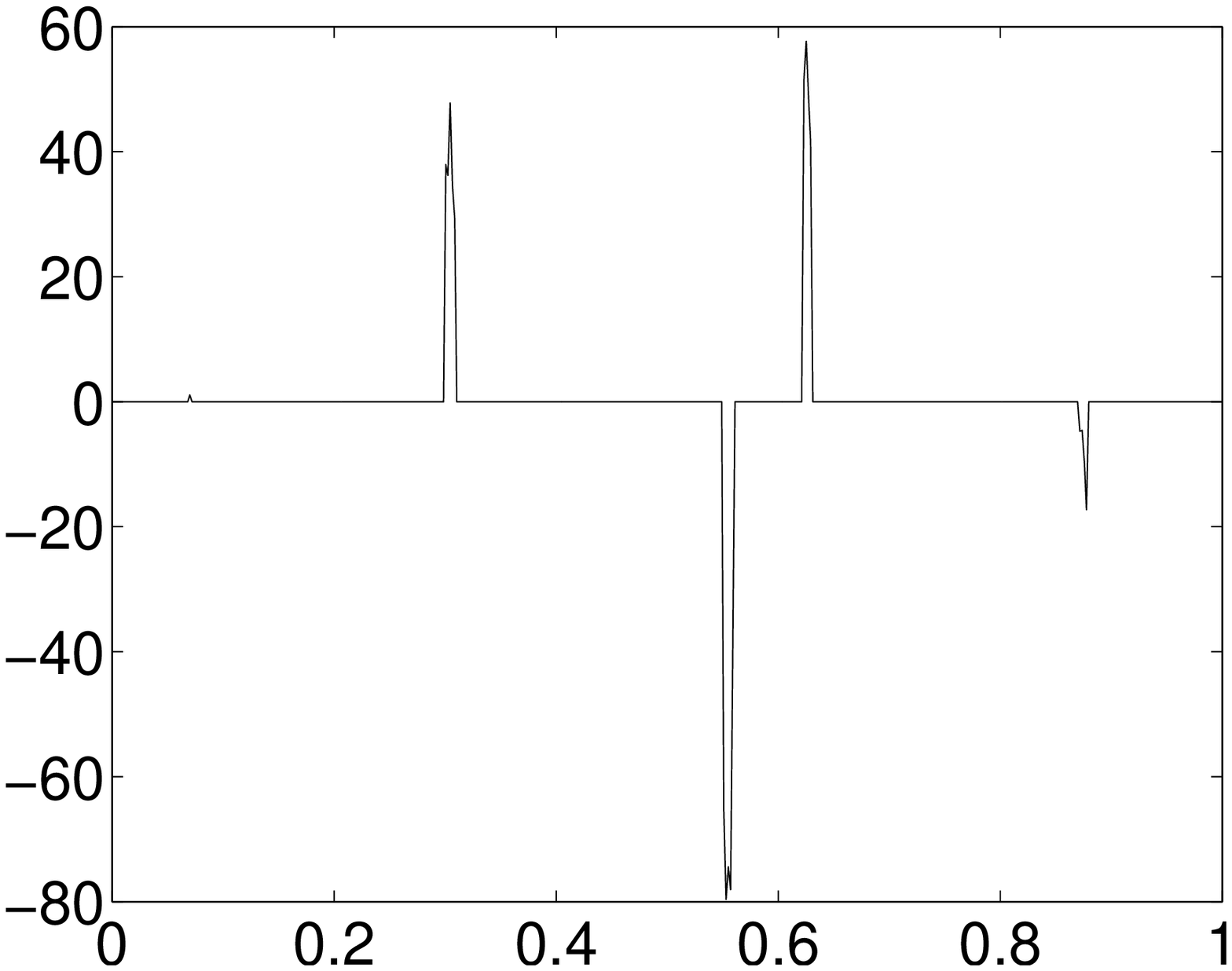}
  \caption{Illustration of data and results of the inverse integration problem.
    Left: the true solution $\bar u$ with $N=500$,
    middle: the noisy data $f$ with 5\% noise,
    right: the reconstruction by $\ell^1$ minimization with $w_k=3 \cdot 10^{-3}$ and $\gamma=5 \cdot 10^5$.
    The solution was obtained with the semismooth Newton method after 11 iterations with a residual norm of $1.7 \cdot 10^{-10}$. }
  \label{fig:illu_inv_integration}
\end{figure}

The true solution $\bar u$ is given by small plateaus and hence the 
data $f^\delta = K\bar u + \delta$ is a noisy function with steep linear ramps.
Figure~\ref{fig:illu_inv_integration} shows our sample data and the result of the $\ell^1$ minimization with the SSN method.
Table~\ref{tab:illu_SSN_inv_integration} shows how the SSN method performed in this specific example.
It can be observed that the residual is not decaying monotonically and it descends slowly in the beginning while it drops significantly in the last step.
Moreover we observed in many examples that the algorithm shows a similar performance for a broad range of starting values $u^0$.
Another important observation is that the performance of the algorithm depends on the value of $\gamma$.
For too small as well as for too large values of $\gamma$ the algorithm does only converge when started very close to the solution.
As a rule of thumb one could take $\gamma$ close to the reciprocal of the smallest singular value of the (in practice unknown) matrix $\KK_{\AA\AA}$ where $\AA$ is the sparsity pattern of the solution.

\begin{table}[htb]
  \centering
  \begin{tabular}{@{}ccc@{}}
 \toprule 
 $n$ & $\Psi(u^n)$ & $\|r^n\|$ \\
 \midrule 
 \verb|  1| & \verb|1.3249e+01| & \verb|6.9764e+05| \\
 \verb|  2| & \verb|1.0461e+01| & \verb|2.1698e+02| \\
 \verb|  3| & \verb|5.3849e+00| & \verb|2.9586e+02| \\
 \verb|  4| & \verb|4.8393e+00| & \verb|8.3922e+02| \\
 \verb|  5| & \verb|4.2488e+00| & \verb|1.9864e+02| \\
 \verb|  6| & \verb|3.0433e+00| & \verb|1.5474e+02| \\
 \verb|  7| & \verb|2.8758e+00| & \verb|3.9127e+01| \\
 \verb|  8| & \verb|2.8237e+00| & \verb|3.5658e+01| \\
 \verb|  9| & \verb|2.7365e+00| & \verb|2.9485e+01| \\
 \verb| 10| & \verb|2.5984e+00| & \verb|7.7932e+00| \\
 \verb| 11| & \verb|2.5518e+00| & \verb|1.7423e-10| \\
 \bottomrule 
\end{tabular}
\bigskip
  \caption{Illustration of the performance of the SSN method for the inverse integration problem. The second column shows the decay of the function value $\Psi$ while the third column shows the norm of the residual. The data is the same as in Figure~\ref{fig:illu_inv_integration}.}
  \label{tab:illu_SSN_inv_integration}
\end{table}

We made experiments to see how the SSN method depend on the noise level
and the regularization parameter.
First, we fixed the regularization sequence and changed the noise level.
Hence, we solved the problem~(\ref{eq:min_prob_inv_int}) for fixed $N=500$,
fixed $w_k = 10^{-5}$ and varied the noise level $\delta$.
Table~\ref{tab:compare_noise_inv_integration} reports the results.
Basically, a higher noise level leads to a smaller number of iterations but longer CPU-time (this is, because the active sets are larger during the iteration).
\begin{table}[htb]
  \centering
  \begin{tabular}{@{}cccc@{}}
 \toprule
 $\|\delta\|$ &  \#iter & CPU-Time (sec) & $\|r\|$ \\
\midrule
 \verb|1.0e+00| &  \verb|  5| & \verb|6.66e-01|  & \verb|6.18e-10|\\
 \verb|1.0e-01| &  \verb|  8| & \verb|7.63e-01|  & \verb|2.19e-09|\\
 \verb|1.0e-02| &  \verb| 12| & \verb|3.98e-01|  & \verb|7.85e-10|\\
 \verb|1.0e-03| &  \verb| 10| & \verb|3.11e-01|  & \verb|1.29e-09|\\
 \verb|1.0e-04| &  \verb| 11| & \verb|3.22e-01|  & \verb|9.85e-10|\\
 \verb|1.0e-05| &  \verb| 11| & \verb|3.18e-01|  & \verb|2.49e-09|\\
\bottomrule
\end{tabular}
\bigskip
  \caption{Behavior of the SSN method for different noise levels with fixed $w_k = 10^{-5}$. The problem under consideration is the inverse integration, the problem size is $N=500$ with $\gamma = 5\cdot 10^5$ throughout.
  The rightmost column shows the residual norm at convergence.}
  \label{tab:compare_noise_inv_integration}
\end{table}
Second we coupled the noise level and the regularizing sequence $w_k$.
Since it is shown in \cite{daubechies2004iteratethresh,lorenz2008reglp}
that a parameter choice $w_k \propto \delta$ provides a regularization
we used $w_k=\delta$.
Hence, we solved the problem~(\ref{eq:min_prob_inv_int}) for fixed $N=500$ and different noise levels $\delta$, see Table~\ref{tab:compare_regparameter_inv_integration} for the results.
Basically, the algorithm behaves similar for different noise levels,
especially the CPU-time is always comparably small.
\begin{table}[htb]
  \centering
  \begin{tabular}{@{}cccc@{}}
 \toprule
 $\|\delta\|$ &  \#iter & CPU-Time (sec) & $\|r\|$ \\
\midrule
 \verb|1.0e-01| &  \verb|  9| & \verb|1.25e-01|  & \verb|6.94e-12|\\
 \verb|1.0e-02| &  \verb| 11| & \verb|2.61e-01|  & \verb|2.29e-11|\\
 \verb|1.0e-03| &  \verb| 11| & \verb|2.96e-01|  & \verb|1.47e-10|\\
 \verb|1.0e-04| &  \verb| 11| & \verb|3.10e-01|  & \verb|2.91e-11|\\
 \verb|1.0e-05| &  \verb| 11| & \verb|3.25e-01|  & \verb|7.63e-11|\\
 \verb|1.0e-06| &  \verb| 12| & \verb|3.50e-01|  & \verb|4.16e-11|\\
\bottomrule
\end{tabular}
\bigskip
  \caption{Behavior of the SSN method for different noise levels. The problem under consideration is the inverse integration, the problem size is $N=500$.
  We chose $w_k = \norm{\delta}$ and $\gamma = 5\cdot 10^5$ throughout.
  The rightmost column shows the residual norm at convergence.}
  \label{tab:compare_regparameter_inv_integration}
\end{table}

\begin{figure}[htb]
  \centering
  \includegraphics[width=.6\textwidth]{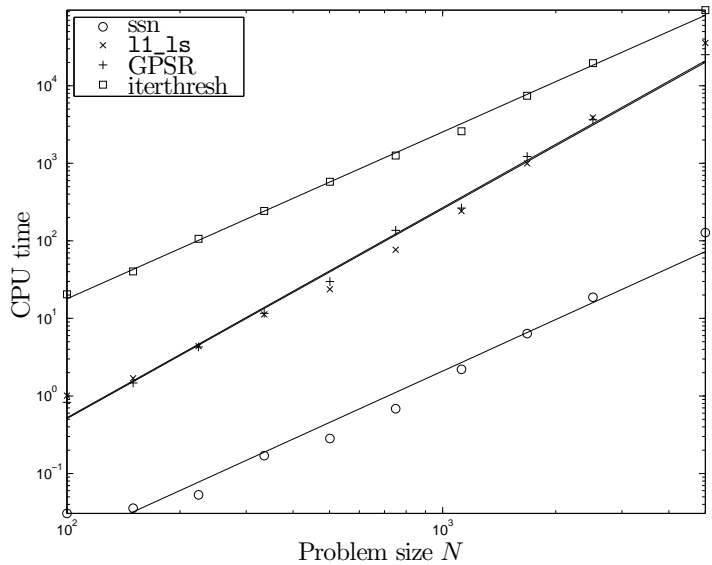}
  \caption{The empirical growth of the computational cost for the different algorithms.}
  \label{fig:compare_methods_inv_integration}
\end{figure}

Moreover, we made a simple experiment to assess how the computational cost grow with the size of the problem.
We considered the inverse intergration problem with problem size $N$ between 100 and 5000.
We kept all parameters, as well as the data and the noise level fixed and only refined the discretization of the problem.
We stopped the algorithms when a required residual tolerance was reached.
Moreover, we checked if the reached functional value was equal for the different methods since the algorithms used different stopping criteria.
Table~\ref{tab:compare_methods_inv_integration} reports CPU times required for the SSN method, for GPSR, \verb!l1_ls! and for the iterative thresholding.
In Figure~\ref{fig:compare_methods_inv_integration} the same data is coded graphically.
When assuming that the computational cost is $\mathcal{O}(N^\beta)$ we found $\beta=2.71$ for GPRS, $\beta = 2.70$ for \verb!l1_ls!, $\beta=2.15$ for iterative thresholding, and $\beta= 2.20$ for SSN.
Moreover, the constant hidden in the $\mathcal{O}$ notation is considerably smaller for the SSN method.
The observed scaling differs from the results reported in~\cite{figueiredo2007gradproj} which may be due to the different structure of the examples.
In~\cite{figueiredo2007gradproj} the example used a matrix which had all singular values either close to one or zero, while in our example the singular values converge to zero.
Hence, it is expected that the empirical scaling of the computational costs differs from problem to problem.

\begin{table}[htb]
  \centering
  \begin{tabular}{@{}ccccc@{}}
 \toprule
 $N$ & SSN & GPRS & \verb!l1_ls! & iterthresh\\
\midrule
 \verb|  100| & \verb|3.06e-02| & \verb|8.29e-01| & \verb|1.01e+00| & \verb|2.04e+01| \\
 \verb|  150| & \verb|3.57e-02| & \verb|1.47e+00| & \verb|1.69e+00| & \verb|4.03e+01| \\
 \verb|  224| & \verb|5.31e-02| & \verb|4.28e+00| & \verb|4.40e+00| & \verb|1.06e+02| \\
 \verb|  335| & \verb|1.70e-01| & \verb|1.16e+01| & \verb|1.12e+01| & \verb|2.43e+02| \\
 \verb|  501| & \verb|2.83e-01| & \verb|2.99e+01| & \verb|2.37e+01| & \verb|5.78e+02| \\
 \verb|  750| & \verb|6.84e-01| & \verb|1.36e+02| & \verb|7.69e+01| & \verb|1.26e+03| \\
 \verb| 1122| & \verb|2.20e+00| & \verb|2.65e+02| & \verb|2.42e+02| & \verb|2.59e+03| \\
 \verb| 1679| & \verb|6.37e+00| & \verb|1.22e+03| & \verb|9.95e+02| & \verb|7.42e+03| \\
 \verb| 2512| & \verb|1.87e+01| & \verb|3.64e+03| & \verb|3.88e+03| & \verb|1.95e+04| \\
 \verb| 5000| & \verb|1.28e+02| & \verb|2.51e+04| & \verb|3.56e+04| & \verb|9.45e+04| \\
\bottomrule
\end{tabular}
\bigskip
  \caption{Comparison of the CPU time in seconds for the different algorithms and different sizes of the problem.
  The problem under consideration  is the inverse integration problem with 5\% noise and regularization parameter $w_k=3 \cdot 10^{-3}$.}
  \label{tab:compare_methods_inv_integration}
\end{table}

\subsection{Deblurring in a Haar basis}
\label{sec:debl-haar-basis}

As a second example of an ill-posed problem we consider a blurring operator
$A:L^2([0,1])\to L^2([0,1])$ given by $Au = k*u$ with the kernel $k(x) = c \, (1+x^2/\lambda^2)^{-1}$ with $\lambda = 0.01$.
We choose $c$ such that $\int k\, dx = 1$
and consider $u$ to be extended periodically to $\R$ in order to evaluate the convolution integral.

In this example we work with a synthesis operator $B:\ell^2 \to L^2([0,1])$ mapping coefficients $(c_k)$ to a function $u=\sum_k c_k \psi_k$ where the $(\psi_k)$ form the orthonormal Haar wavelet basis~\cite{maass1997wavelet}.
Hence, the operator under consideration $K = AB$ is a blurring after a Haar wavelet synthesis, see~\cite{chambolle1998shrink,daubechies2004iteratethresh} for discussions of $\ell^1$ penalty terms in combination with wavelet expansions.

\begin{figure}[htb]
  \centering
  \includegraphics[width=.4\textwidth]{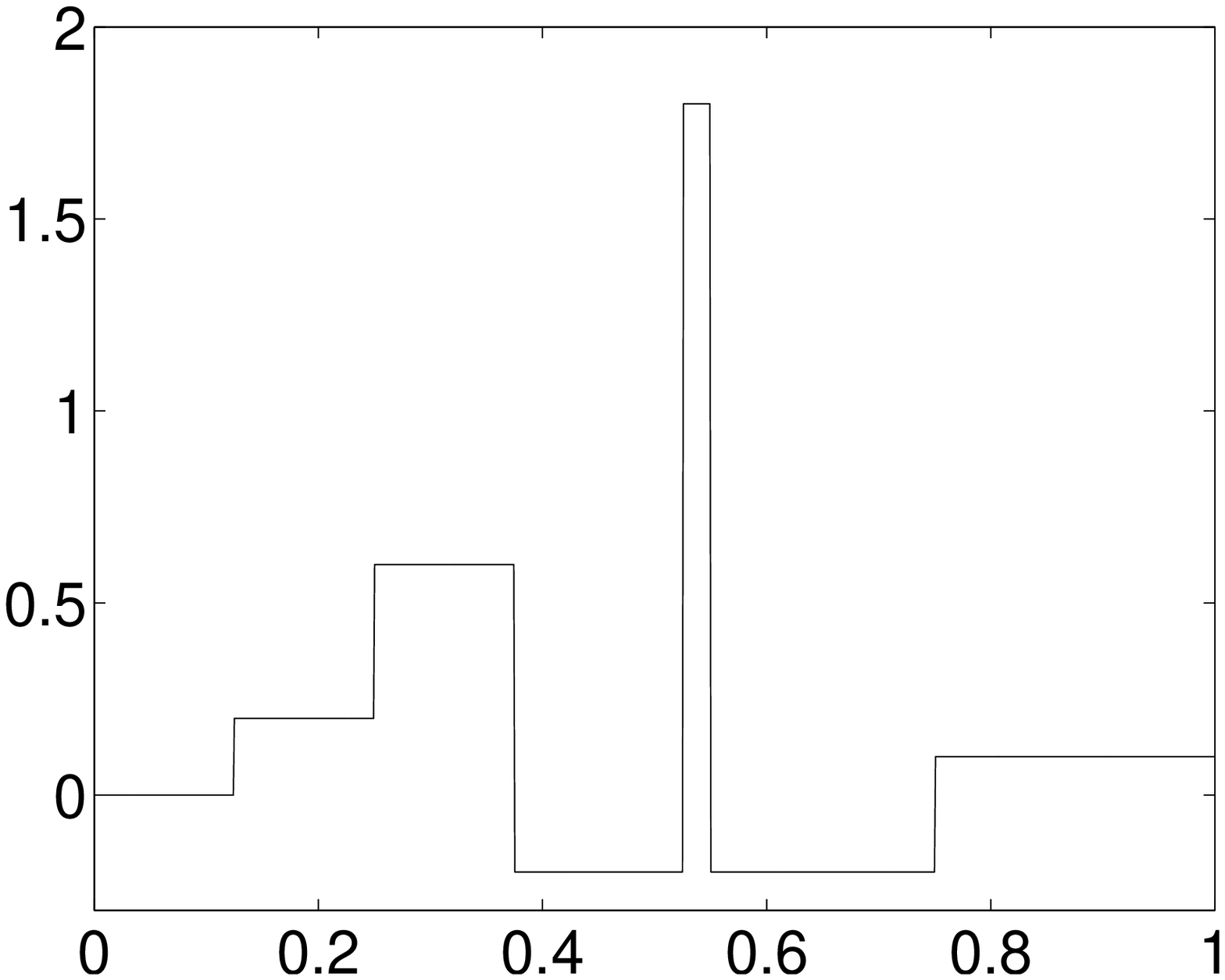}\ 
  \includegraphics[width=.4\textwidth]{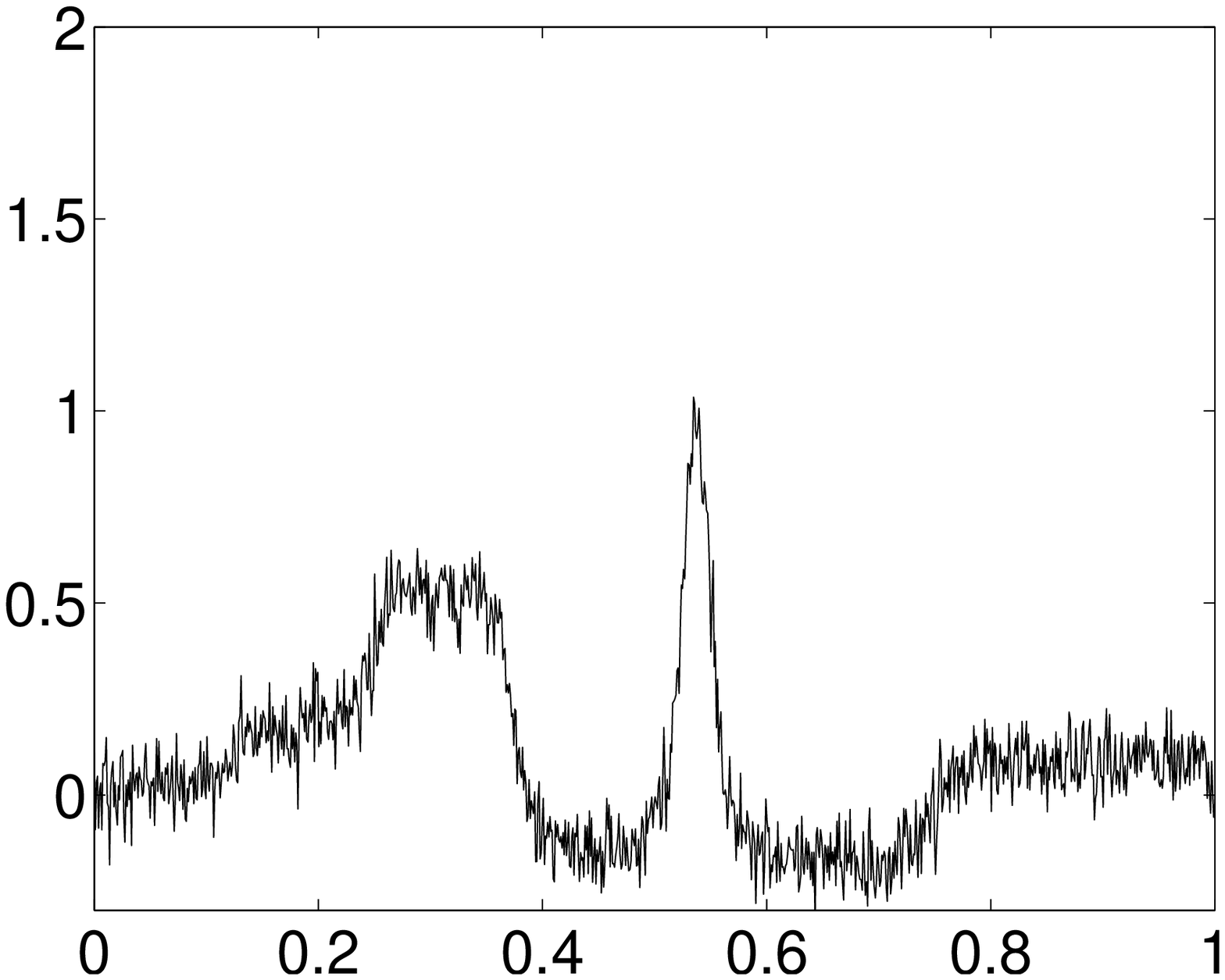}\\
  \includegraphics[width=.4\textwidth]{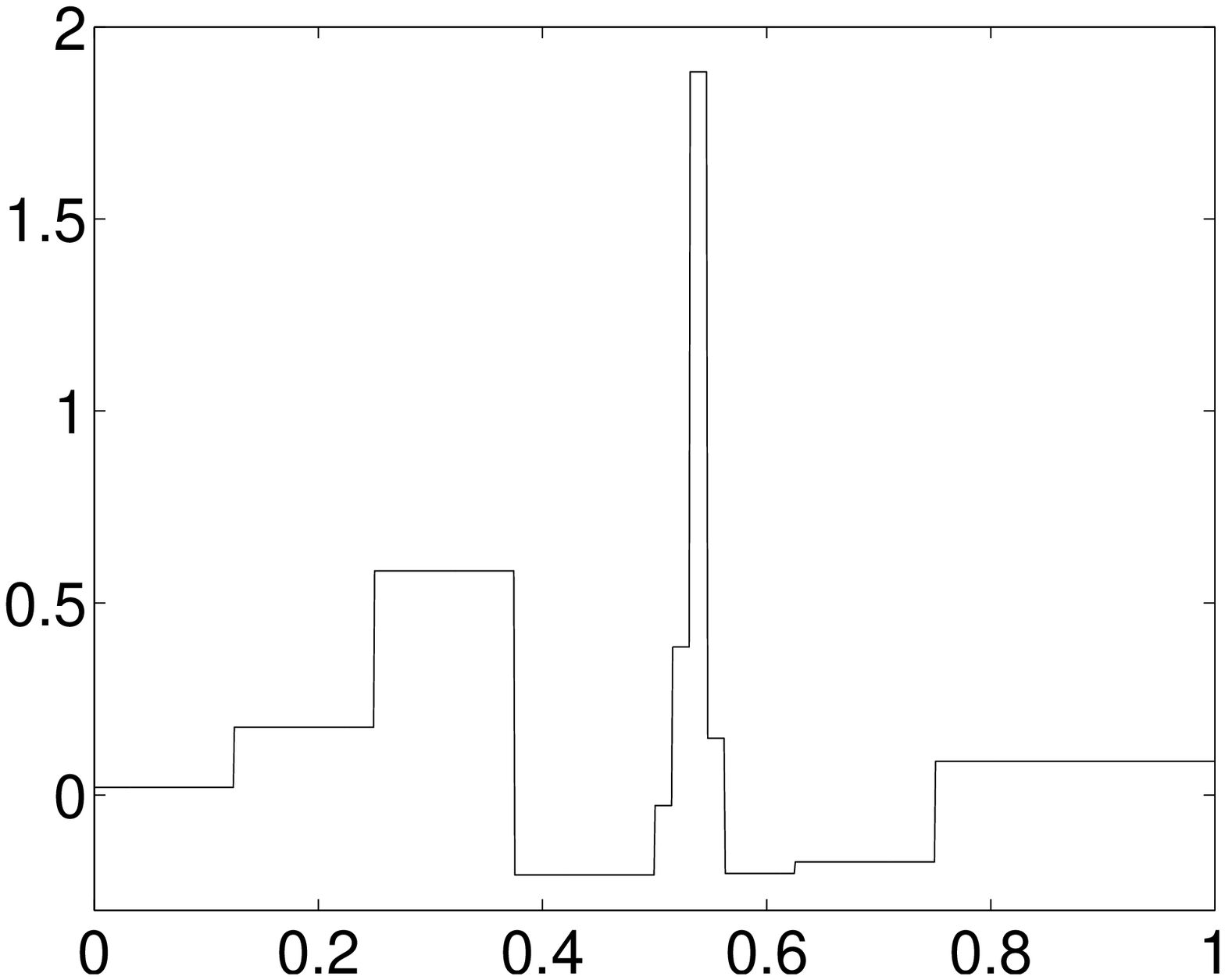}\ 
  \includegraphics[width=.4\textwidth]{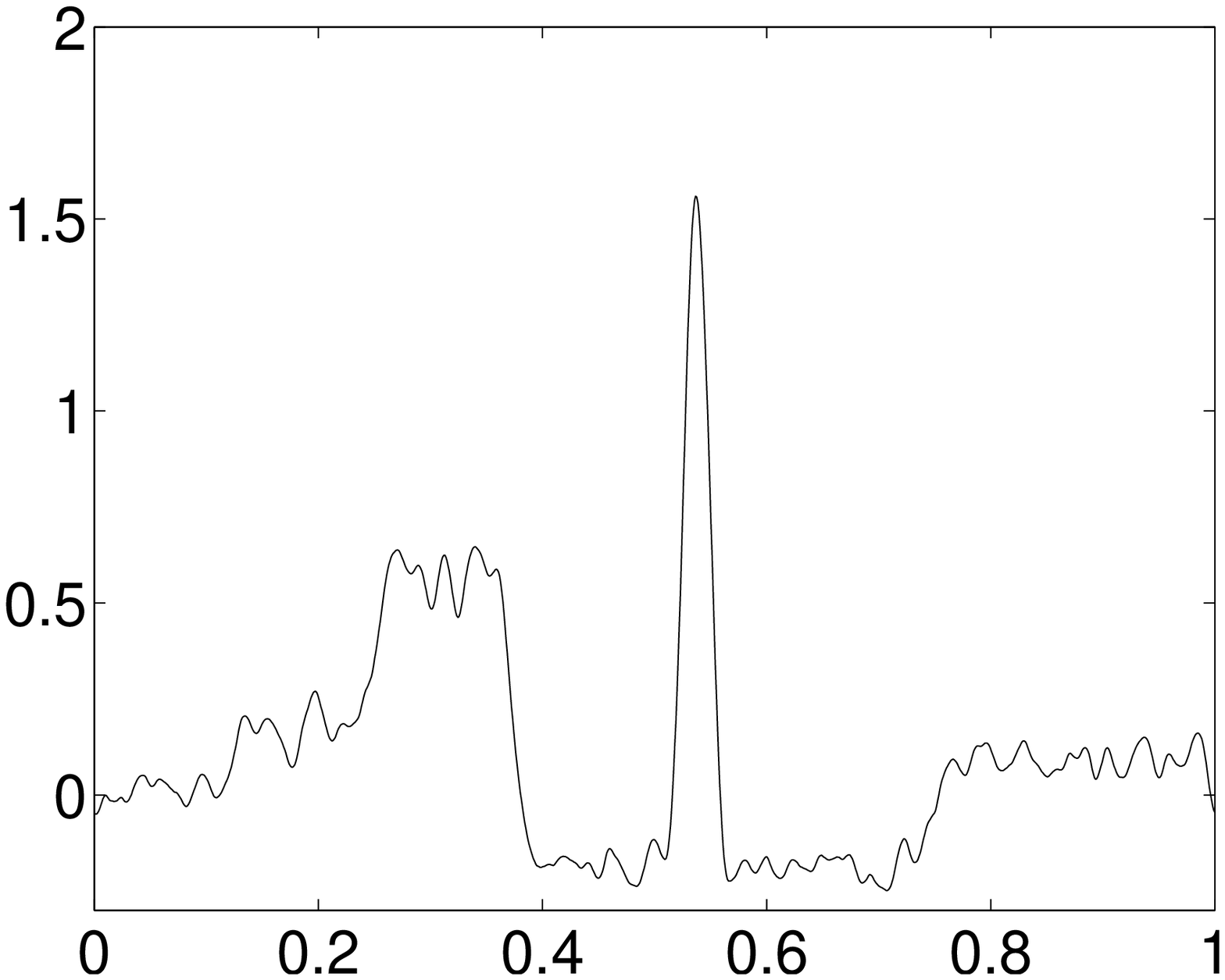}
  \caption{The results of $\ell^1$ and $\ell^2$ (classical Tikhonov) regularization of deblurring in a Haar basis. Upper left: the true solution $u$, upper right: the given data $f$, lower left: the reconstruction by $\ell^1$ minimization with $w_k = 0.12$ and $\gamma=5\cdot 10^6$, lower right: the reconstruction by $\ell^2$ minimization with $w_k = 0.05$.}
  \label{fig:smoothing_example}
\end{figure}

We start with a given function $u$ which is piecewise constant.
The data $f$ is computed as $f = Au + \text{noise}$ such that we have 25\% relative error, i.e.~ $\norm{f-Au}/\norm{f} = 0.25$.
The Haar coefficients of $u$ have been reconstructed by minimizing~(\ref{eq:Main_problem}).
As an illustration of $\ell^1$ penalties in contrast to classical $\ell^2$ regularization we also show the results of the minimization of
\[
  \frac{1}{2} \norm{Kc - f}^2 +
  \sum_{k=1}^\infty w_k |c_k|^2 \ .
\]
Figure~\ref{fig:smoothing_example} and Table~\ref{tab:illu_smoothing} show the results of both $\ell^1$ and the above $\ell^2$ regularization where we discretized the problem to 1024 Haar wavelets.
The parameters $w_k$ are independent of $k$ and have been tuned by hand to produce optimal results.
Since the original data is quite sparse in the Haar wavelet basis, the $\ell^1$-reconstruction leads to much better results, as expected from the model.
It also turned out that the algorithm is robust with respect to different initial values $u^0$.
We tested several initial values (starting at zero, at $K^* f$ or at a random position) and the observed convergence behavior was very similar in all cases. 

The SSN method converged in six iterations and in 0.3 seconds (for comparison: the GPSR method takes 0.5 seconds and \verb!l1_ls! converged in 5.3 seconds). 

\begin{table}[htb]
  \centering
  \begin{tabular}{@{}ccc@{}}
 \toprule 
 $n$ & $\Psi(u^n)$ & $\|r^n\|$ \\
 \midrule 
 \verb|  1| & \verb|3.3920e+001| & \verb|2.9676e+006| \\
 \verb|  2| & \verb|1.3905e+002| & \verb|3.2499e+004| \\
 \verb|  3| & \verb|1.3326e+001| & \verb|6.2647e+005| \\
 \verb|  4| & \verb|7.9347e+000| & \verb|1.7517e+004| \\
 \verb|  5| & \verb|6.0006e+000| & \verb|8.0510e-002| \\
 \verb|  6| & \verb|5.9823e+000| & \verb|1.5424e-009| \\
 \bottomrule 
\end{tabular}
\bigskip
  \caption{Illustration of the performance of the SSN method for deblurring in a Haar basis. The second column shows the decay of the function value $\Psi$ while the third column shows the norm of the residual. The data is the same as in Figure~\ref{fig:smoothing_example}.}
  \label{tab:illu_smoothing}
\end{table}

\subsection{Compressive Sampling}
\label{sec:compressive-sampling}
In our last example we illustrate the applicability of the SSN method to the decoding problem in compressive sampling alias compressed sensing (CS).
In CS one aims at reconstructing a signal from very few linear measurements, see\cite{candes2006compressivesampling,donoho2006compressedsensing} for an introduction to CS.
A popular way of decoding a signal from data $f$ which was measured by the observation operator $K$ is to minimize a functional of type~(\ref{eq:Main_problem}), see \cite{candes2005decoding}.
Our example on compressive sampling is taken from~\cite{figueiredo2007gradproj}.
We obtain an observation operator $K\in\R^{K\times N}$ by first filling it with independent samples of a standard Gaussian distribution and then orthonormalizing the rows.
Hence, the operator is not injective but it possesses the so-called restricted isometry property (see~\cite{baraniuk2007riprandommatrices}) which means that all submatrices consisting of a small number of columns have singular values close to one.
Especially, submatrices made of a small number of columns are injective.
Hence, the SSN method works as long as the active sets are small enough.

\begin{figure}[htb]
  \centering
  \includegraphics[width=.4\textwidth]{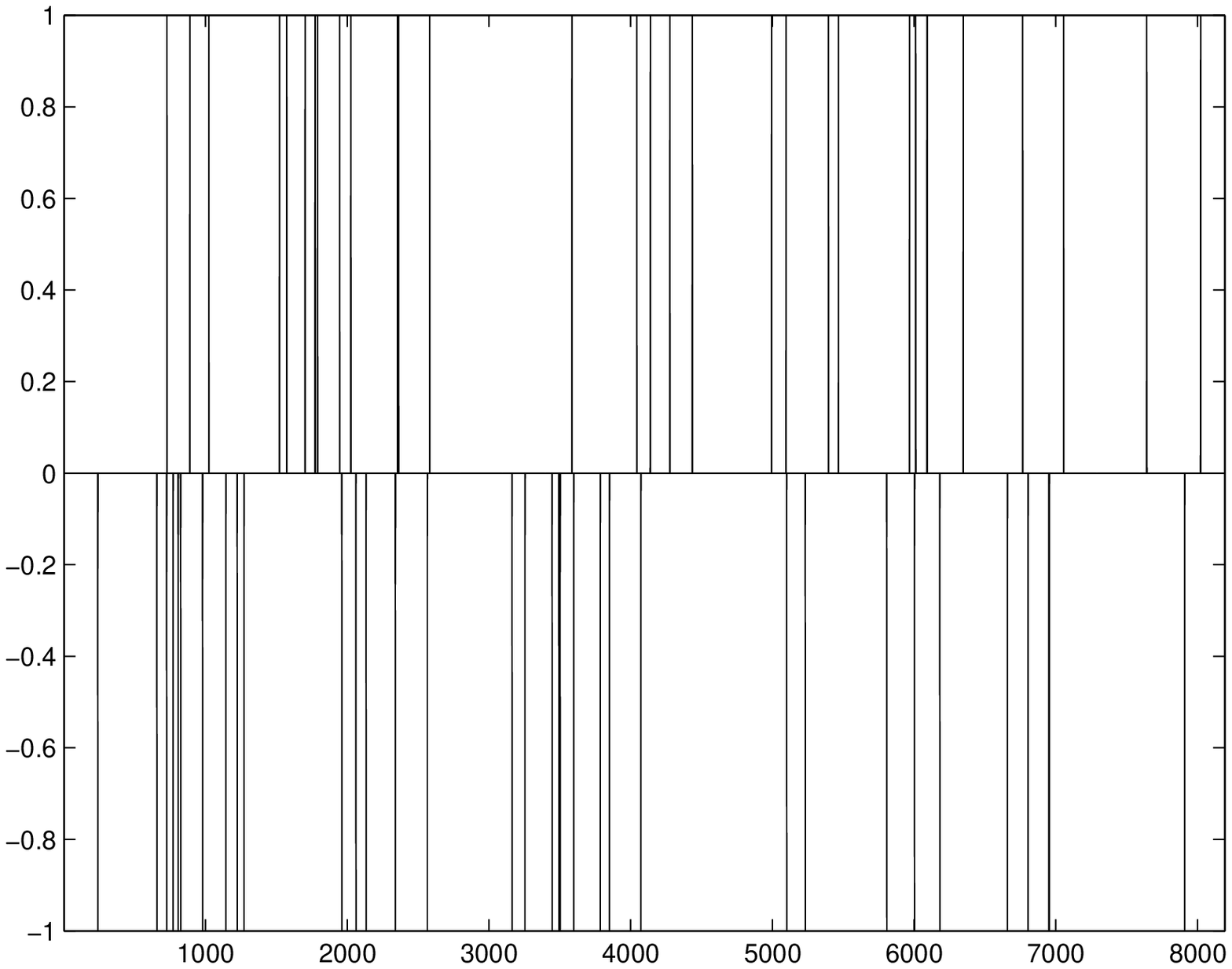}
  \includegraphics[width=.4\textwidth]{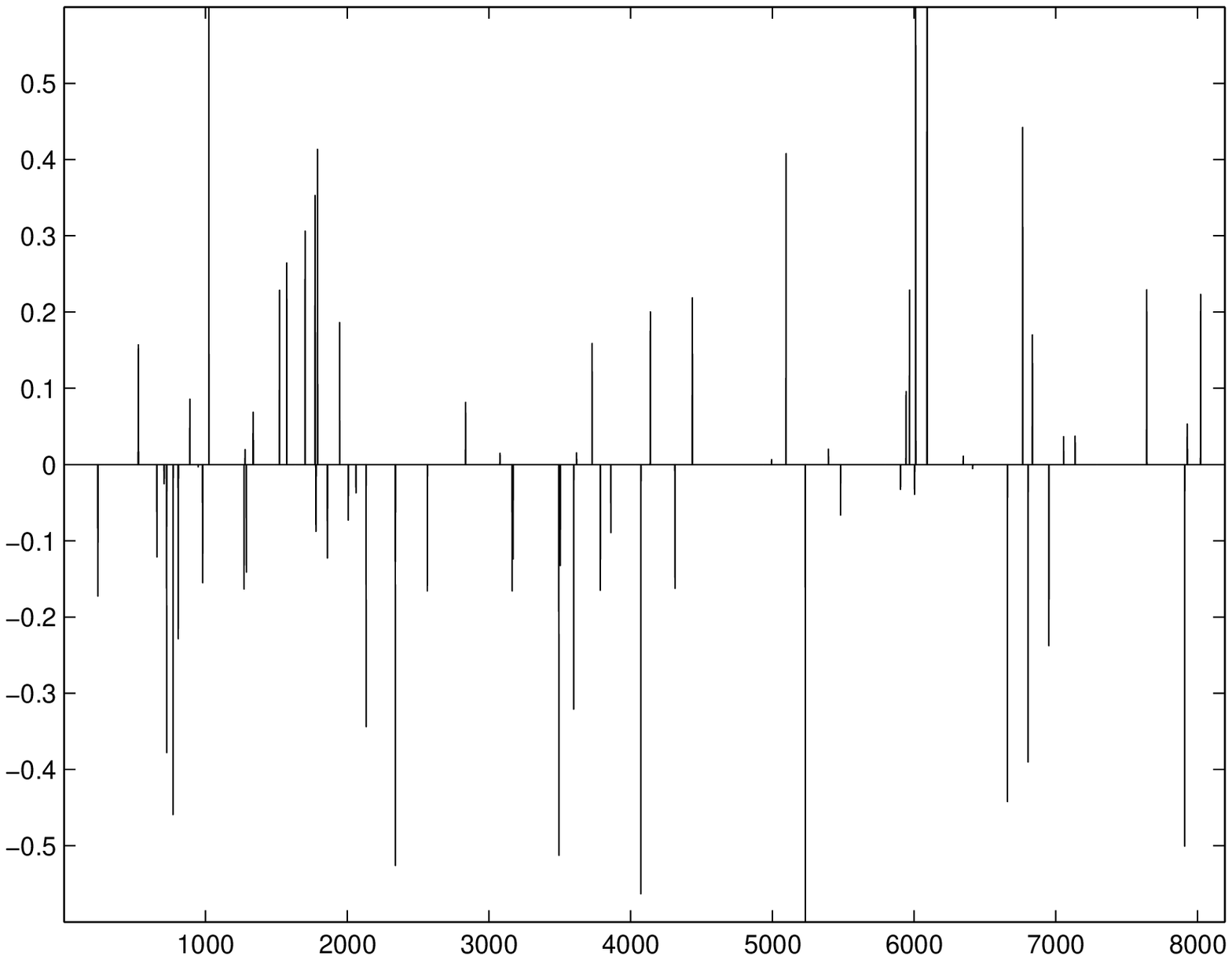}
  \caption{Illustration of data and results of the CS example.
    Left: the original signal $u$ with $n=8192$,
    right: the reconstruction by $\ell^1$ minimization with $w_k=0.05$ and $\gamma=5\cdot 10^4$.
    The solution was obtained with the semismooth Newton method after 6 iterations with a residual norm of approximately $1\cdot 10^{-11}$. }
  \label{fig:illu_SSN_cs}
\end{figure}

In this example we chose $N=8192$, $K=512$, and the signal $u$ contained $64$ randomly placed $\pm 1$ spikes.
The observation $f$ was generated by $f = Ku + \text{noise}$ such that we have 5\% relative error.
The minimization of~(\ref{eq:Main_problem}) with $w=0.05$ was done with the SSN method with parameter $\gamma = 5 \cdot 10^{4}$.
The SSN method converged in approximately 1.2 seconds in six iterations and the active sets stayed very small during the iteration, see Figure~\ref{fig:illu_SSN_cs} and Table~\ref{tab:illu_SSN_cs}.
Hence, the SSN method is a promising candidate for the decoding problem in CS.

\begin{table}[htb]
  \centering
  \begin{tabular}{@{}ccccc@{}}
\toprule
 $n$ & $\Psi(u^n)$ & $\|r^n\|$ & $|\AA|$ & $\text{cond}(\KK_{\AA\AA})$ \\\midrule
 \verb$ 1$ & \verb$2.0774e+00$ & \verb$1.2254e+04$ & \verb|252| & \verb| 44.52|\\
 \verb$ 2$ & \verb$7.1752e+00$ & \verb$2.2715e+03$ & \verb|148| & \verb| 12.77|\\
 \verb$ 3$ & \verb$2.7379e+00$ & \verb$4.6644e+02$ & \verb| 90| & \verb|  6.31|\\
 \verb$ 4$ & \verb$1.9997e+00$ & \verb$1.4674e+02$ & \verb| 67| & \verb|  4.75|\\
 \verb$ 5$ & \verb$1.8386e+00$ & \verb$3.9728e+01$ & \verb| 67| & \verb|  4.60|\\
 \verb$ 6$ & \verb$1.8361e+00$ & \verb$9.9652e-12$ & \verb| 67| & \verb|      |\\
\bottomrule
\end{tabular}
\bigskip
  \caption{Illustration of the performance of the SSN method for CS. The second column shows the decay of the function value $\Psi$ while the third column shows the norm of the residual. The forth and fifth column show the size of the active set and the condition of the matrix $\KK_{\AA\AA}$ which has to be inverted in the Newton step. The data is the same as in Figure~\ref{fig:illu_SSN_cs}.}
  \label{tab:illu_SSN_cs}
\end{table}


\section{Conclusion}
\label{sec:Conclusion}
We have shown that the semismooth Newton method applied to Tikhonov functionals with sparsity constraints is a fast algorithm which is easy to implement as an active set method.
Each iteration involves the solution of a system of linear equations on the active coefficients only.
Our numerical experiments show that these systems stay reasonably small during the iteration and are also very well conditioned.
In addition, the experiments indicate that the SSN method compares favorably with existing state-of-the-art methods when applied to ill-posed problems.
While we investigated only the local convergence behavior, the numerical experiments indicate that our method is robust with respect to the initial value of the iteration.
However, the convergence is slow as long as the iterates are far from the minimizer and it gets faster when the solution is approached.
The global convergence properties are not yet explained by our theory and need further investigation.
Another direction for further research is globalization of the method e.g., by the use of an appropriate merit function, and line search or trust region methods.


\begin{appendix}
\section*{Appendix}
\label{sec:Appendix}
We define for $u \in \ell^2$ and $h \in \HH$
\begin{equation*}
  F(u) = \sum_{k=1}^\infty w_k \, |u_k|, \qquad
  G(h) = \frac{1}{2} \norm{h-f}^2_\HH
\end{equation*}
and calculate their conjugate (polar) functions, see \cite[Ch.~I.4]{ekelandtemam1976convex}.
We have
\begin{align*}
  F^*(p) & = \sup_{u \in \ell^2} \big( \langle p, u \rangle - F(u) \big) 
  = \sup_{u} \big( \langle p, u \rangle - \sum_{k=1}^\infty w_k \, |u_k| \big) \\
  & = \sup_{u} \big( \sum_{k=1}^\infty (p_k - w_k \text{ sign } u_k) \, u_k \big) 
  =
  \begin{cases}
    0, & \text{if } |p_k| \le w_k \text{ for all } k \\
    \infty & \text{otherwise}.
  \end{cases}
\end{align*}
For $G$, we obtain
\begin{equation*}
  G^*(p) = \sup_{h \in \HH} \big( \langle p,h \rangle - G(h) \big) 
  = \sup_h \big( \langle p,h \rangle - \frac{1}{2} \norm{h-f}^2_\HH \big) 
  = \frac{1}{2} \norm{p}^2_\HH + \langle p, f \rangle,
\end{equation*}
since the supremum is attained at $h = p + f$.


\end{appendix}

\section*{References}

\bibliographystyle{plain}
\bibliography{SSN_Sparsity_Constraints}

\end{document}